\documentclass{amsart}[12pt]

\theoremstyle{plain}
\begingroup
\newtheorem{thm}{Theorem}[section]
\newtheorem{lem}[thm]{Lemma}
\newtheorem{cor}[thm]{Corollary}

\endgroup

\theoremstyle{definition}
\newtheorem{defn}[thm]{Definition}

\theoremstyle{remark}
\newtheorem{rem}[thm]{Remark}

\begin{document}

\def\headline #1{\medskip \centerline{\bf #1} \medskip}

\def\K {\operatorname{K}\nolimits}
\def\H {\operatorname{H}\nolimits}
\def\height {\operatorname{height}\nolimits}
\def\Ker {\operatorname{Ker}\nolimits}
\def\Im {\operatorname{Im}\nolimits}
\def\Image {\operatorname{Image}\nolimits}
\def\initial {\operatorname{in}\nolimits}
\def\Gin{\operatorname{Gin}\nolimits}
\def\Spec{\operatorname{Spec}\nolimits}
\def\D{\operatorname{D}\nolimits}
\def\V{\operatorname{V}\nolimits}
\def\S{\operatorname{S}\nolimits}
\def\E{\operatorname{E}\nolimits}
\def\rank{\operatorname{rank}\nolimits}
\def\HT{\operatorname{HT}\nolimits}

\title
[The Hilbert series of algebras of Veronese type]
{The Hilbert series of algebras of Veronese type}
\author{Mordechai Katzman}
\address{
Department of Pure Mathematics,\\ University of Sheffield,\\Hicks Building,\\ Sheffield S3 7RH,\\ United Kingdom\\
e-mail: M.Katzman@sheffield.ac.uk } \maketitle
\begin{abstract}
This paper gives an explicit formula for the Hilbert series of algebras of Veronese type.
\end{abstract}

\section{Introduction}

In this paper we describe the Hilbert series of algebras of Veronese type:

\begin{defn}
Fix a positive integer $d$ and a sequence of integers
$\mathbf{a}=(a_1,\dots,a_n)$ such that
$1\le a_1\le \dots \le a_n\le d$ and $\sum_{i=1}^n a_i >d$.
Let ${\mathcal V}(\mathbf{a};d)$ be the $k$-subalgebra of $k[x_1,\dots,x_n]$
generated by all monomials $x_1^{\alpha_1} \dots x_n^{\alpha_n}$ with
$\sum_{i=1}^n \alpha_i=d$ and $\alpha_i\le a_i$ for all $1\le i \le n$.

We shall also denote by $\mathcal S$ all subsets $S$ of $\{1,\dots,n\}$
with $\sum_{i\in S} a_i < d$; for any $S\subset \{1,\dots,n\}$ we define
$\Sigma S$ to be $\sum_{i\in S} a_i$.
\end{defn}

Note that
${\mathcal V}(d,d,\dots,d; d)$ is the classical Veronese algebra
while
${\mathcal V}(1,1,\dots,1; d)$ is the monomial algebra
associated with the $d$th hypersimplex.

For the purpose of computing Hilbert series and $a$-invariants we will
use a normalized grading on these algebras so that the degree of their
generators equals one.

These monomial algebras have recently atracted considerable interest;
E.~DeNegri,  T.~Hibi (\cite{DH}) have recently classified those which
are Gornestein and B.~Sturmfels
(\cite{Stu}) described Gr\"obner bases arising from presentations of
these algebras.

It is known that algebras of Veronese type are normal and in \cite{DH}
the authors classified all such algebras which are Gorenstein.
Also, in \cite{BVV} the authors described the canonical modules
and $a$-invariants of  ${\mathcal V}(1,1,\dots,1; d)$.
The aim of this section is to extend and complement these results
by producing an explicit formula for the $h$-vectors
of all algebras of Veronese type. Additionally, the explicit
formulas provide a very efficient way for computing these Hilbert
series.

The results in this paper are part of an ongoing project to understand certain monomial algebras.
Since these results have been obtained several articles have referred to them in a preliminary form (e.g., \cite{NOH});
it is hoped that this paper will make the results more widely available.

\section{The Hilbert series of algebras of Veronese type}

We start by describing the Hilbert function of algebras of Veronese type:

\begin{thm}\label{Vhfun}
The Hilbert function of  ${\mathcal V}(\mathbf{a};d)$ is given by
$$H(i)=\sum_{S\in {\mathcal S}} (-1)^{|S|}
\binom{i(d-\Sigma S)-|S|+n-1}{n-1}.$$
\end{thm}

\begin{proof}

By Lemma 2.1 in \cite{DH} $H(i)$ is the number of sequences
$(\alpha_1,\dots,\alpha_n)$ satisfying $\sum_{j=1}^n \alpha_j = id$ and
$1\le \alpha_j \le i a_i$. This number is readily seen to be the coefficient of
$t^{id}$ in
$$\prod_{i=1}^n (1+t+\dots+t^{i a_i})=
\frac{\prod_{i=1}^n 1-t^{i a_i +1}}{(1-t)^n}=$$
$$\left( \sum_{S\subset \{1,\dots,n\}} (-1)^{|S|} t^{i \Sigma S+|S|}\right)
\left( \sum_{r=0}^\infty \binom{r+n-1}{n-1} t^r \right)$$
and the coefficient of
$t^{id}$ in this expression is
$$\sum_{S\in {\mathcal S}} (-1)^{|S|} \binom{i(d-\Sigma S)-|S|+n-1}{n-1}$$
\end{proof}

\begin{cor}
\label{Ehrhart}
The Hilbert function $H_{nd}(i)$ of ${\mathcal V}(1,1,\dots,1; d)$ is given by
$$H_{nd}(i)=\sum_{s=0}^{d-1} (-1)^s \binom{n}{s}
\binom{i(d-s)-s+n-1}{n-1}.$$
\end{cor}

\begin{rem}
Since ${\mathcal V}(1,1,\dots,1; d)$ are normal, the previous Corollary also gives the
Ehrhart polynomial of the $d$th hypersimplex. This
was also proved for $d=2$ in
Chapter 9  of \cite{Stu} using Gr\"obner bases techniques.

T.~Hibi has also noted that one could obtain the Hilbert functions of
algebras of Veronese type
using their initial ideals discussed in
the last chapter of \cite{Stu}.
\end{rem}

To produce the $h$-vectors of
${\mathcal V}(\mathbf{a};d)$ we need to understand
the generating functions of $\binom{n+id-1}{n-1}$, and these will
be described using the following:

\begin{defn}
For any positive integers $n$ and $d$ define the numbers $A^{n,d}_i$ by
$$\left( 1+T+\dots+T^{d-1} \right)^n = \sum_{i\ge 0} A^{n,d}_i T^i .$$
\end{defn}

\begin{thm}\label{genfun}
$$(1-t)^n \sum_{i=0}^\infty \binom{n+id-1}{n-1} t^i =
\sum_{j\ge 0} A^{n,d}_{jd} t^j$$
\end{thm}

\begin{proof}

Let $\Xi$ be the set of (complex) $d$th roots of $1$; we have
\begin{gather*}
\frac{1}{d} \sum_{\xi\in \Xi} \frac{1}{(1-\xi t^{1/d})^n}= \\
\frac{1}{d} \sum_{\xi\in \Xi}
 \sum_{j=0}^\infty \binom{n+j-1}{n-1} (\xi t^{1/d})^j=\\
\sum_{j=0}^\infty \binom{n+jd-1}{n-1} t^j ,
\end{gather*}
the last equality following from the fact that
$$\sum_{\xi\in \Xi} \xi^j=
\left\{
\begin{array}{l l}
d & d | j \\
0 & \text{otherwise}.
\end{array}
\right.
$$
Multiplying both sides by $(1-t)^n=(1-(\xi t^{1/d})^d)^n$ we obtain
\begin{gather*}
(1-t)^n \sum_{j=0}^\infty \binom{n+jd-1}{n-1} t^j=\\
\frac{1}{d} \sum_{\xi\in \Xi}
 \frac{\left(1-(\xi t^{1/d})^d\right)^n} {\left( 1-\xi t^{1/d} \right)^n}= \\
\frac{1}{d} \sum_{\xi\in \Xi}
 \left( 1+\xi t^{1/d}+ \dots + (\xi t^{1/d})^{d-1} \right)^n.
\end{gather*}
The coefficient of $(t^{1/d})^s$ in this expression is
$$\frac{1}{d} \sum_{\xi\in \Xi} \xi^s A^{n,d}_s=
\left\{
\begin{array}{l l}
A^{n,d}_s & d | s\\
0 & \text{otherwise,}
\end{array}
\right.
$$
thus only integer powers of $t$ appear in the sum and the coefficient of
$t^s$ is $A^{n,d}_{sd}$.
\end{proof}

The previous Theorem provides an elementary proof for
Theorem 2.4(a) in \cite{DH}:

\begin{cor}
The $h$-vector of ${\mathcal V}(d,d,\dots,d; d)$
is symmetric if and only if $d$ divides $n$,
therefore, ${\mathcal V}(d,d,\dots,d; d)$
is Gorenstein if and only if $d$ divides $n$.
The multiplicity is $d^{n-1}$ and the $a$-invariant equals to
$\displaystyle -\left\lceil \frac{n}{d} \right\rceil$.
\end{cor}

\begin{proof}
Note that the $h$-vector
is given by the coefficients of
$$(1-t)^n \sum_{i=0}^\infty \binom{n+id-1}{n-1} t^i=
\sum_{i\ge 0} A^{n,d}_{id} T^i $$
and that
$ \sum_{i\ge 0} A^{n,d}_{id} T^i$
is a unimodal and reciprocal
polynomial of degree
$n(d-1)$ (cf. section 3.5 in \cite{A}.)

If $d$ divides $n$ then
$\sum_{i\ge 0} A^{n,d}_{id} T^i$ is obviously reciprocal. On the other
hand, if $d$ does not divide $n$ then $d$ does not divide $n(d-1)$ either
and in this case the coefficient of the highest power of $T$ in
$\sum_{i\ge 0} A^{n,d}_{id} T^i$ will be greater than $1=A^{n,d}_0$.

We compute the $a$-invariant as the degree of the Hilbert series as
a rational function: the highest power of $t$ occurring in the numerator
is
$$\left\lfloor\frac{n(d-1)}{d}\right\rfloor=
n+\left\lfloor\frac{-n}{d}\right\rfloor=
n-\left\lceil\frac{n}{d}\right\rceil .
$$

The multiplicity will be computed in Corollary \ref{mult} below.
\end{proof}

\begin{lem}
For any non-negative integers $n,a$ and $b$
let
$$P^n_{a,b}(t)=\sum_{i=0}^{\infty} \binom{ai-b+n-1}{n-1} t^i .$$
For any positive integer $b$ we have
$P^n_{a,b}=\sum_{j=0}^b (-1)^j \binom{b}{j} P^{n-j}_{a,0}$.
\end{lem}

\begin{proof}
Since
$$\binom{ai-b+n-1}{n-1}=
\binom{ai-(b-1)+n-1}{n-1}-\binom{ai-(b-1)+n-2}{n-2}$$
we have a recursion relation
$P^n_{a,b}=P^n_{a,b-1}-P^{n-1}_{a,b-1}$
and the result follows.
\end{proof}

Combining this with Theorem \ref{Vhfun} we obtain an explicit
expression for the $h$-vectors of ${\mathcal V}(\mathbf{a}; d)$:

\begin{thm} \label{Vhvec}
\begin{gather*}
(1-t)^n \sum_{i=0}^\infty
\sum_{S\in {\mathcal S}} (-1)^{|S|}
\binom{i(d-\Sigma S)-|S|+n-1}{n-1} t^i= \\
\sum_{S\in {\mathcal{S}}} (-1)^{|S|}
\sum_{j=0}^{|S|} (-1)^j \binom{|S|}{j}
(1-t)^j \sum_{l\ge 0} A^{n-j,d-\Sigma S}_{l(d-\Sigma S)} t^l .
\end{gather*}
\end{thm}

\begin{proof}
For any $S\subset \{1,\dots,n\}$ we have
\begin{gather*}
(1-t)^n  (-1)^{|S|}
\sum_{i=0}^\infty  \binom{i(d-\Sigma S)-|S|+n-1}{n-1} t^i=\\
(1-t)^n   (-1)^{|S|} P^n_{d-\Sigma S,|S|}=\\
(-1)^{|S|} \sum_{j=0}^{|S|}
(-1)^j \binom{|S|}{j} (1-t)^j ((1-t)^{n-j} P^{n-j}_{d-\Sigma S,0})
\end{gather*}
and by Theorem \ref{genfun} this equals
$$ (-1)^{|S|}
\sum_{j=0}^{|S|} (-1)^j \binom{|S|}{j}
(1-t)^j \sum_{l\ge 0} A^{n-j,d-\Sigma S}_{l(d-\Sigma S)} t^l .$$
\end{proof}

\begin{cor}
The Hilbert series of ${\mathcal V}(1,1,\dots,1; d)$ is
$$
(1-t)^{-n}
\sum_{s=0}^{d-1} (-1)^s \binom{n}{s} \sum_{j=0}^s (-1)^j \binom{s}{j}
(1-t)^j \sum_{l\ge 0} A^{n-j,d-s}_{l(d-s)} t^l .$$
For $d=2$ this reduces to
$$ (1-t)^{-n} \left( \sum_{l\ge 0} \binom{n}{2l} t^l  -nt \right) .$$
\end{cor}

\begin{proof}
The first statement follows easily from the previous Theorem. To prove
the second statement note that $A^{n,2}_j=\binom{n}{j}$ and that
$A^{n,1}_j=0$ unless $j=0$, in which case we have $A^{n,1}_0=1$.
\end{proof}

Even though the expressions for the $h$-vectors for
${\mathcal V}({\mathbf a};d)$ look forbidding, it is possible
to extract useful information from them. We shall now proceed to
compute the multiplicities and find a (sharp)
upper bound for their $a$-invariants.

\begin{lem}
$$\sum_{l\ge 0} A^{n,d}_{ld}=d^{n-1} .$$
\end{lem}

\begin{proof}
Let $\Xi$ be the set of (complex) $d$th roots of one. We have
\begin{gather*}
\sum_{l\ge 0} A^{n,d}_{ld} t^l =\\
\frac{1}{d} \sum_{\xi\in \Xi} \sum_{l\ge 0} A^{n,d}_{l}
\left(\xi t^{1/d}\right)^l=\\
\frac{1}{d}  \sum_{\xi\in \Xi} \left( 1+ \xi t^{1/d} +
\left(\xi t^{1/d}\right)^2+\dots+\left(\xi t^{1/d}\right)^{d-1} \right)^n .
\end{gather*}
Setting $t=1$ the last sum vanishes unless $\xi=1$ and, therefore,
the sum equals
$d^{n-1}$.
\end{proof}

\begin{cor}\label{mult}
The multiplicity of  ${\mathcal V}(\mathbf{a}; d)$ is
$$\sum_{S\in {\mathcal S}} (-1)^{|S|} (d-\Sigma S)^{n-1} .$$
\end{cor}

\begin{proof}
Substituting $t=1$ in the numerator of the Hilbert series
obtained in Theorem \ref{Vhvec} we obtain
$$\sum_{S\in {\mathcal S}} (-1)^{|S|}
\sum_{l\ge 0} A^{n,d-\Sigma S}_{l(d-\Sigma S)}  $$
which by the previous Lemma equals
$$\sum_{S\in {\mathcal S}} (-1)^{|S|} (d-\Sigma S)^{n-1}$$
\end{proof}

\begin{cor}
If $n\ge d$ the $a$-invariant of
${\mathcal V}(\mathbf{a}; d)$ is at most
$\displaystyle -\left\lceil \frac{n}{d} \right\rceil$.
\end{cor}

\begin{proof}
The $a$-invariant is the degree of the Hilbert series
of ${\mathcal V}(\mathbf{a}; d)$ as a rational function. Note
that the highest degree of $t$ occurring in a summand of
$$
\sum_{S\in {\mathcal{S}}} (-1)^{|S|}
\sum_{j=0}^{|S|} (-1)^j \binom{|S|}{j}
(1-t)^j \sum_{l\ge 0} A^{n-j,d-\Sigma S}_{l(d-\Sigma S)} t^l $$
is
\begin{gather*}
\max_{S\in {\mathcal S}} \max_{0\le j\le |S|}
j+\left\lfloor \frac{(n-j)(d-\Sigma S-1)}{d-\Sigma S}\right\rfloor=\\
j+n-j- \min_{S\in {\mathcal S}} \min_{0\le j\le |S|}
\left\lceil \frac{n-j}{d-\Sigma S}\right\rceil=\\
n- \min_{S\in {\mathcal S}} \left\lceil \frac{n-|S|}{d-\Sigma S}\right\rceil
\le \\
n- \min_{S\in {\mathcal S}} \left\lceil \frac{n-|S|}{d- |S|}\right\rceil=
n -\left\lceil \frac{n}{d} \right\rceil
\end{gather*}
the last equality holding for $n\ge d$.
\end{proof}


\end{document}